\documentclass[]{amsart}

\usepackage{enumitem}
\usepackage{newclude}
\usepackage{multicol}
\usepackage[labelsep=period]{caption}
\usepackage{tabu}
\usepackage[table]{xcolor}
\usepackage{tikz}
\usetikzlibrary{arrows,automata,positioning}
\usepackage{rotating}
\usepackage{graphicx}
\usepackage{eurosym}
\usepackage{amsfonts}
\usepackage{amsmath}
\usepackage{amssymb}
\usepackage{stmaryrd}
\usepackage{wasysym}
\usepackage{amsthm}
\usepackage[margin=1in]{geometry}
\usepackage[hang,flushmargin,symbol*]{footmisc}
\usepackage{color}
\definecolor{darkblue}{rgb}{0, 0, .6}
\definecolor{grey}{rgb}{.7, .7, .7}
\usepackage[breaklinks]{hyperref}
\hypersetup{
	colorlinks=true,
	linkcolor=darkblue,
	anchorcolor=darkblue,
	citecolor=darkblue,
	pagecolor=darkblue,
	urlcolor=darkblue,
	pdftitle={},
	pdfauthor={}
}

\theoremstyle{definition}
\newtheorem{theorem}{Theorem}[section]

\newtheorem{conjecture}[theorem]{Conjecture}

\newtheorem{definition}[theorem]{Definition}

\newtheorem{proposition}[theorem]{Proposition}

\newsavebox{\savepar}

\begin{document}

\title{Prime Vertex Labelings of Families of Unicyclic Graphs}

\author[Diefenderfer]{Nathan Diefenderfer}
\author[Hastings]{Michael Hastings}
\author[Heath]{Levi N.~Heath}
\author[Prawzinsky]{Hannah Prawzinsky}
\author[Preston]{Briahna Preston}
\author[White]{Emily White}
\author[Whittemore]{Alyssa Whittemore}

\address{Department of Mathematics and Statistics, Northern Arizona University, Flagstaff, AZ 86011}
\thanks{This research was supported by the National Science Foundation grant \#DMS-1148695 through the Center for Undergraduate Research (CURM)}

\date{\today}

\begin{abstract}
A simple $n$-vertex graph has a prime vertex labeling if the vertices can be injectively labeled with the integers $1, 2, 3,\ldots, n$ such that adjacent vertices have relatively prime labels. We will present previously unknown prime vertex labelings for new families of graphs, all of which are special cases of Seoud and Youssef's conjecture that all unicyclic graphs have a prime labeling.
\end{abstract}

\maketitle

\tikzstyle{vert} = [circle, draw, fill=grey,inner sep=0pt, minimum size=4.5mm]
\tikzstyle{b} = [draw,very thick,blue,-]
\tikzstyle{r} = [draw, very thick, red,-]
\tikzstyle{g} = [draw, very thick, green, -]
\tikzstyle{blk} = [draw, very thick, black, -]

\section{Introduction}\label{sec:intro}

Applications of combinatorial graphs can be found everywhere in life, from communication networks to possible moves in a board game. This paper will focus on graph labeling, which is the process of assigning labels to either the vertices, edges, or both, following some predetermined set of rules. More specifically, we study a particular type of labeling of the vertices of a graph, called a prime vertex labeling, where the labels of adjacent vertices are required to be relatively prime. The research presented in this paper was inspired by a conjecture stated in 1999 by Seoud and Youssef~\cite{Seoud1999}, namely:
\begin{quote}
\begin{center}
\emph{All unicyclic graphs have a prime vertex labeling.}
\end{center}
\end{quote}
A unicyclic graph is a graph containing exactly one cycle as a subgraph. Instead of attempting to prove the conjecture outright, which we anticipate would require heavy-duty linear algebra, we focused our attention on finding prime labelings for specific families of unicyclic graphs.

The paper will proceed as follows. In Section~\ref{sec:terminology}, some basics of graph theory will be discussed. Then we will define a prime vertex labeling and state some known results. In Sections~\ref{sec:hairy} and \ref{sec:ternary}, we will present several new families with prime vertex labelings. These new families will consist of graphs having exactly one cycle together with either pendants or pendants with ternary trees attached to each cycle vertex. Finally, in Section~\ref{sec:conclusion}, some conjectures and potential future work will be described.
\section{Graph Theory Terminology}\label{sec:terminology}
This section will provide an overview of the definitions and terminology that will be used throughout the rest of the article.

First, a \emph{graph} $G$ is a set of vertices, $V(G)$, together with a set of edges, $E(G)$, connecting some subset, possibly empty, of the vertices. If $u,v \in V(G)$ are connected by an edge, we say $u$ and $v$ are \emph{adjacent} and the corresponding edge is denoted $uv$ or $vu$. We will restrict our attention to \emph{simple graphs}, which are graphs that do not contain multiple edges between pairs of vertices or have edges that connect a vertex to itself (called a loop). For the remainder of this paper, all graphs are assumed to be simple.

A graph $H$ whose vertex set and edge set are subsets of the vertex set and edge set of a given graph $G$ is a \emph{subgraph} of $G$. We say that the \emph{degree} of a vertex $u$ is the number of edges having $u$ as an endpoint. A graph is \emph{connected} if it does not consist of two or more disjoint ``pieces."

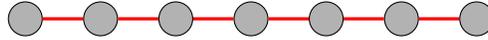
\begin{figure}[!ht]
\begin{center}
\begin{tikzpicture} {scale=1,auto}
\node (a) at (-3,0) [vert] {};
\node (b) at (-2,0) [vert] {};
\node (c) at (-1,0) [vert] {};
\node (d) at (0,0) [vert] {};
\node (e) at (1,0) [vert] {};
\node (f) at (2,0) [vert] {};
\node (g) at (3,0) [vert] {};
\path[r] (a) to (b);
\path[r] (b) to (c);
\path[r] (c) to (d);
\path[r] (d) to (e);
\path[r] (e) to (f);
\path[r] (f) to (g);
\end{tikzpicture}
\end{center}
\caption{The path $P_{7}$.}\label{fig:P7}
\end{figure}

\begin{figure}[!tht]
\begin{center}
\begin{tikzpicture}[scale=1,auto]
\node (a) at (90:2) [vert] {};
\node (b) at (120:2) [vert] {};
\node (c) at (150:2) [vert] {};
\node (d) at (180:2) [vert] {};
\node (e) at (210:2) [vert] {};
\node (f) at (240:2) [vert] {};
\node (g) at (270:2) [vert] {};
\node (h) at (300:2) [vert] {};
\node (i) at (330:2) [vert] {};
\node (j) at (360:2) [vert] {};
\node (k) at (30:2) [vert] {};
\node (l) at (60:2) [vert] {};
\path[r] (a) to (b);
\path[r] (b) to (c);
\path[r] (c) to (d);
\path[r] (d) to (e);
\path[r] (e) to (f);
\path[r] (f) to (g);
\path[r] (g) to (h);
\path[r] (h) to (i);
\path[r] (i) to (j);
\path[r] (j) to (k);
\path[r] (k) to (l);
\path[r] (l) to (a);
\end{tikzpicture}
\end{center}
\caption{The cycle $C_{12}$.}\label{fig:C12}
\end{figure}
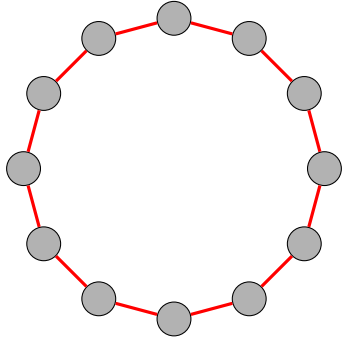

Next, we define a few important families of graphs.  An \emph{$n$-path} (or simply \emph{path}), denoted $P_n$, is the connected graph consisting of two vertices of degree $1$ and $n-2$ vertices of degree $2$. The graph in Figure~\ref{fig:P7} depicts the path $P_7$.  An \emph{$n$-cycle} (or simply \emph{cycle}), denoted $C_n$, is the connected graph consisting of $n$ vertices each of degree $2$. The graph $C_{12}$ is shown in Figure~\ref{fig:C12}. Note that $C_{n}$ always has $n$ vertices and $n$ edges.

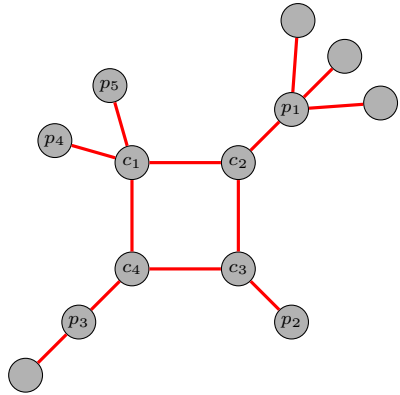
\begin{figure}[!tht]
\begin{center}
\begin{tikzpicture}[scale=1,auto]
\node (1) at (45:1) [vert] {\scriptsize $c_2$};
\node (2) at (135:1) [vert] {\scriptsize $c_1$};
\node (3) at (225:1) [vert] {\scriptsize $c_4$};
\node (4) at (315:1) [vert] {\scriptsize $c_3$};
\node (5) at (45:2) [vert] {\scriptsize $p_1$};
\node (6) at (45:3) [vert] {\scriptsize};
\node (7) at (315:2) [vert] {\scriptsize $p_2$};
\node (8) at (225:2) [vert] {\scriptsize $p_3$};
\node (9) at (225:3) [vert] {\scriptsize};
\node (10) at (150:2) [vert] {\scriptsize $p_4$};
\node (11) at (120:2) [vert] {\scriptsize $p_5$};
\node (12) at (60:3) [vert] {\scriptsize};
\node (13) at (30:3) [vert] {\scriptsize};
\path [r] (1) to (2);
\path [r] (2) to (3);
\path [r] (3) to (4);
\path [r] (4) to (1);
\path [r] (1) to (5);
\path [r] (5) to (6);
\path [r] (4) to (7);
\path [r] (3) to (8);
\path [r] (8) to (9);
\path [r] (2) to (10);
\path [r] (2) to (11);
\path [r] (5) to (12);
\path [r] (5) to (13);
\end{tikzpicture}
\end{center}
\caption{Example of a unicyclic graph consisting of five pendants.}\label{fig:unicyclic pendants}
\end{figure}

As mentioned in Section~\ref{sec:intro}, a major focus of this paper will be \emph{unicyclic graphs}, which have a unique subgraph isomorphic to a cycle. Every vertex lying on the cycle of a unicyclic graph will be referred to as a \emph{cycle vertex}.  In a unicyclic graph, a \emph{pendant} is a path on two vertices with exactly one vertex being a cycle vertex.  The non-cycle vertex of a pendant is called a \emph{pendant vertex}. For example, the graph shown in Figure~\ref{fig:unicyclic pendants} is a unicyclic graph with five pendants. In this case, the vertices labeled by $c_1, c_2, c_3$, and $c_4$ are cycle vertices while the vertices labeled by $p_1, p_2, p_3, p_4$, and $p_5$ are pendant vertices.

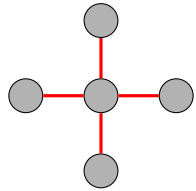
\begin{figure}[!tht]
\begin{center}
\begin{tikzpicture} {scale=1,auto}
\node (a) at (0,0) [vert] {};
\node (b) at (-1,0) [vert] {};
\node (c) at (1,0) [vert] {};
\node (d) at (0,1) [vert] {};
\node (e) at (0,-1) [vert] {};
\path[r] (a) to (b);
\path[r] (a) to (c);
\path[r] (a) to (d);
\path[r] (a) to (e);
\end{tikzpicture}
\end{center}
\caption{The star $S_{4}$.}\label{fig:S4}
\end{figure}

An \emph{$n$-star} (or simply \emph{star}), denoted $S_n$, is the graph consisting of one vertex of degree $n$ and $n$ vertices of degree $1$. Note that $S_n$ consists of $n+1$ vertices and $n$ edges. The star $S_4$ is shown in Figure~\ref{fig:S4}.

A \emph{tree} is a graph having no subgraph isomorphic to a cycle. One defining characteristic of trees is that there exists exactly one ``trail" of edges between every pair of vertices. Paths and stars are examples of trees.

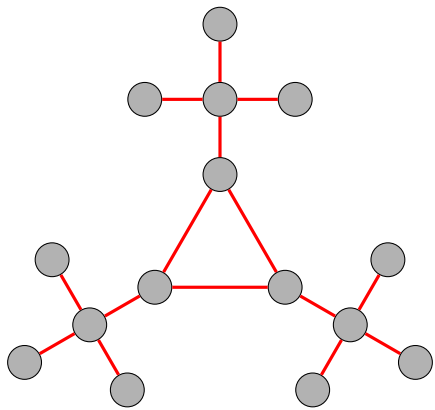
\begin{figure}[!ht]
\begin{center}
\begin{tikzpicture}[scale=1,auto]
\node (a) at (90:1) [vert] {};
\node (b) at (90:2) [vert] {};
\node (c) at (90:3) [vert] {};
\node (d) at (-1,2) [vert] {};
\node (e) at (1,2) [vert] {};
\path[r] (a) to (b);
\path[r] (b) to (c);
\path[r] (b) to (e);
\path[r] (b) to (d);
\begin{scope}[rotate=120]
\node (a') at (90:1) [vert] {};
\node (b') at (90:2) [vert] {};
\node (c') at (90:3) [vert] {};
\node (d') at (-1,2) [vert] {};
\node (e') at (1,2) [vert] {};
\path[r] (a') to (b');
\path[r] (b') to (c');
\path[r] (b') to (e');
\path[r] (b') to (d');
\end{scope}
\begin{scope}[rotate=-120]
\node (a'') at (90:1) [vert] {};
\node (b'') at (90:2) [vert] {};
\node (c'') at (90:3) [vert] {};
\node (d'') at (-1,2) [vert] {};
\node (e'') at (1,2) [vert] {};
\path[r] (a'') to (b'');
\path[r] (b'') to (c'');
\path[r] (b'') to (e'');
\path[r] (b'') to (d'');
\end{scope}
\path[r] (a) to (a');
\path[r] (a') to (a'');
\path[r] (a'') to (a);
\end{tikzpicture}
\end{center}
\caption{An example of the unicyclic graph $C_3\star P_2 \star S_3$.}\label{fig:unicyclic example}
\end{figure}

Most graphs in this paper will result from ``selectively gluing" copies of trees to the cycle vertices of an $n$-cycle. For example, the graph in Figure~\ref{fig:unicyclic example} is a unicyclic graph that results from attaching a copy of the path $P_2$ to each cycle vertex of $C_3$ followed by attaching copies of the star $S_3$ at the vertex of degree 3 to each of the pendant vertices. This particular graph will be denoted by $C_3\star P_2 \star S_3$.  Note that our $\star$ notation is particular to this paper and is not a construction typically found in the literature.

\section{Graph Labeling}\label{sec:labeling}

There are a plethora of ways to label graphs, and their basic premises are similar: assign numbers to the edges or vertices that follow specified rules. Graph labelings have many different applications, including cryptography, wireless networking, radar, and even radio astronomy. The focus of this paper is on prime vertex labelings.

Recall that two integers $a$ and $b$ are said to be \emph{relatively prime} if their greatest common factor is $1$, denoted $(a,b)=1$.

\begin{definition}\label{def:prime labeling}
A simple graph with $n$ vertices is said to have a \emph{prime vertex labeling} (or simply a \emph{prime labeling}) if there is an injection $f:V\rightarrow \{1,2, \ldots n\}$ such that for each edge $uv \in E(G)$, $(f(u),f(v))=1$. For brevity, if a graph has a prime vertex labeling, we will say that the graph is \emph{coprime}.
\end{definition}

\begin{figure}[h!]
\begin{center}
\begin{tikzpicture}[scale=1,auto]
\node (a) at (0,0) [vert] {1};
\node (b) at (1,0) [vert] {2};
\node (c) at (2,0) [vert] {3};
\node (d) at (3,0) [vert] {4};
\node (e) at (0,-1) [vert] {6};
\node (f) at (1,-1) [vert] {7};
\node (g) at (2,-1) [vert] {8};
\node (h) at (3,-1) [vert] {5};
\path[r] (a) -- (b) -- (c) -- (d);
\path[r] (e) -- (f) -- (g) -- (h);
\path[r] (a) to (e);
\path[r] (b) to (f);
\path[r] (c) to (g);
\path[r] (d) to (h);
\path[r] (a) to (f);
\end{tikzpicture}
\end{center}
\caption{An example of a prime labeled graph.}\label{fig:labeling1}
\end{figure}
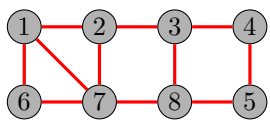

The graph in Figure~\ref{fig:labeling1} depicts one possible prime labeling.

When attempting to find or identify prime vertex labelings, the following basic facts from number theory are useful:
\begin{itemize}
\item All pairs of consecutive integers are relatively prime;
\item All pairs of consecutive odd integers are relatively prime;
\item A common divisor of two numbers is also a divisor of their difference;
\item The integer $1$ is relatively prime to all integers.
\end{itemize}

One motivation for our research is the following conjecture first made by Seoud and Youssef in~\cite{Seoud1999}.

\begin{conjecture}\label{conj:unicyclic}
All unicyclic graphs are coprime.
\end{conjecture}
\section{Known Prime Vertex Labelings}\label{sec:known}

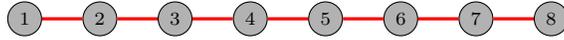
\begin{figure}[!ht]
\begin{center}
\begin{tikzpicture}[scale=1,auto]
\node (a) at (0,0) [vert] {\scriptsize $1$};
\node (b) at (1,0) [vert] {\scriptsize $2$};
\node (c) at (2,0) [vert] {\scriptsize $3$};
\node (d) at (3,0) [vert] {\scriptsize $4$};
\node (e) at (4,0) [vert] {\scriptsize $5$};
\node (f) at (5,0) [vert] {\scriptsize $6$};
\node (g) at (6,0) [vert] {\scriptsize $7$};
\node (h) at (7,0) [vert] {\scriptsize $8$};
\path [r] (a) -- (b) -- (c) -- (d) -- (e) -- (f) -- (g) -- (h);
\end{tikzpicture}
\end{center}
\caption{The path $P_8$ with a prime labeling.}\label{fig:path}
\end{figure}

Prime vertex labelings are known to exist for several infinite families of graphs. We will discuss a few notable labelings that are relevant to the graphs we study in this paper. It has been shown that every path $P_n$ is coprime~\cite{Gallian2014}. In particular, one can use the obvious linear ordering, as seen in Figure~\ref{fig:path}.

\begin{figure}
\begin{center}
\begin{tikzpicture}[scale=.75,auto]
\node (a) at (90:2) [vert] {\scriptsize $1$};
\node (b) at (120:2) [vert] {\scriptsize $12$};
\node (c) at (150:2) [vert] {\scriptsize $11$};
\node (d) at (180:2) [vert] {\scriptsize $10$};
\node (e) at (210:2) [vert] {\scriptsize $9$};
\node (f) at (240:2) [vert] {\scriptsize $8$};
\node (g) at (270:2) [vert] {\scriptsize $7$};
\node (h) at (300:2) [vert] {\scriptsize $6$};
\node (i) at (330:2) [vert] {\scriptsize $5$};
\node (j) at (360:2) [vert] {\scriptsize $4$};
\node (k) at (30:2) [vert] {\scriptsize $3$};
\node (l) at (60:2) [vert] {\scriptsize $2$};
\path[r] (a) to (b);
\path[r] (b) to (c);
\path[r] (c) to (d);
\path[r] (d) to (e);
\path[r] (e) to (f);
\path[r] (f) to (g);
\path[r] (g) to (h);
\path[r] (h) to (i);
\path[r] (i) to (j);
\path[r] (j) to (k);
\path[r] (k) to (l);
\path[r] (l) to (a);
\end{tikzpicture}
\end{center}
\caption{The cycle $C_{12}$ with a prime labeling.}\label{fig:cycle}
\end{figure}
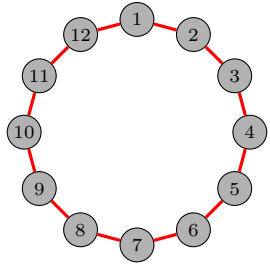

Similarly, every cycle $C_n$ can be labeled using the obvious linear ordering so that the vertices labeled $1$ and $n$ are adjacent~\cite{Gallian2014}. For example, see the labeling given in Figure~\ref{fig:cycle}.

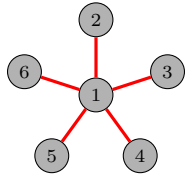
\begin{figure}
\begin{center}
\begin{tikzpicture}[scale=1,auto]
\node (a) at (0:0) [vert] {\scriptsize $1$};
\node (b) at (18:1) [vert] {\scriptsize $3$};
\node (c) at (90:1) [vert] {\scriptsize $2$};
\node (d) at (162:1) [vert] {\scriptsize $6$};
\node (e) at (234:1) [vert] {\scriptsize $5$};
\node (f) at (-54:1) [vert] {\scriptsize $4$};
\path[r] (a) to (b);
\path[r] (a) to (c);
\path[r] (a) to (d);
\path[r] (a) to (e);
\path[r] (a) to (f);
\end{tikzpicture}
\end{center}
\caption{The star $S_5$ with a prime labeling.}\label{fig:star}
\end{figure}

A prime labeling for every star $S_n$ can be constructed by labeling the unique vertex of degree $n$ with 1 and the remaining vertices with the integers 2 through $n$ in any order~\cite{Gallian2014}, as seen in Figure~\ref{fig:star}.

The infinite family of graphs that consist of a cycle with a path of length $m$ attached to each cycle vertex, denoted $C_n\star P_m$, are also coprime~\cite[Theorem 2.5]{Seoud1999}.

Additionally, the graph constructed from first attaching a pendant to every cycle vertex of $C_n$, then attaching a complete binary tree (i.e., a directed rooted tree with every internal vertex having two children) at each pendant vertex has a prime labeling~\cite[Theorem 2.6]{Seoud1999}. These graphs are the inspiration for our investigation of cycles with pendants having complete ternary trees attached to the pendant vertices that will be presented in Section~\ref{sec:ternary}.

Other examples of infinite families of graphs that are known to be coprime include complete graphs $K_n$ if and only if $n\leq3$~\cite{Gallian2014}, wheels $W_n$ if and only if $n$ is even ~\cite{Gallian2014}, all helms $H_n$, and all books $B_n$~\cite[Theorem 2.3]{Seoud1999}. Consult Gallian's dynamic survey~\cite{Gallian2014} for a comprehensive listing of the families of graphs that are known to have or known not to have prime vertex labelings. In~\cite{Seoud2012}, Seoud et al.~provide necessary and sufficient conditions for a graph to be coprime, but we will not elaborate on that here.
\section{Hairy Cycles}\label{sec:hairy}

In this section, the first set of new results will be shown. All of the graphs are constructed by attaching pendants to the vertices of a cycle. These graphs are called hairy cycles.

\begin{definition}
For all $m, n \in \mathbb{N}$ with $n\geq 3$, an \emph{$m$-hairy $n$-cycle}, denoted $C_n \star S_m$, is the cycle $C_{n}$ with $m$ pendants attached to each cycle vertex.
\end{definition}

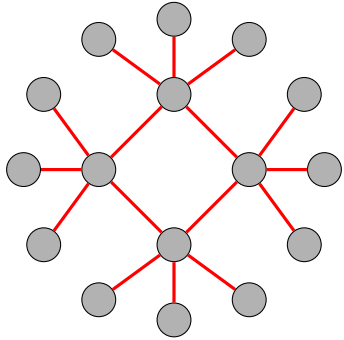
\begin{figure}[!ht]
\begin{center}
\begin{tikzpicture}[scale=1,auto]
\node (1) at (90:1) [vert] {\scriptsize };
\node (2) at (120:2) [vert] {\scriptsize };
\node (3) at (90:2) [vert] {\scriptsize };
\node (4) at (60:2) [vert] {\scriptsize };
\begin{scope}[rotate=90]
\node (15) at (90:1) [vert] {\scriptsize };
\node (13) at (120:2) [vert] {\scriptsize };
\node (14) at (90:2) [vert] {\scriptsize };
\node (16) at (60:2) [vert] {\scriptsize };
\end{scope}
\begin{scope}[rotate=-90]
\node (7) at (90:1) [vert] {\scriptsize };
\node (5) at (120:2) [vert] {\scriptsize };
\node (6) at (90:2) [vert] {\scriptsize };
\node (8) at (60:2) [vert] {\scriptsize };
\end{scope}
\begin{scope} [rotate=180]
\node (11) at (90:1) [vert] {\scriptsize };
\node (9) at (120:2) [vert] {\scriptsize };
\node (10) at (90:2) [vert] {\scriptsize };
\node (12) at (60:2) [vert] {\scriptsize };
\end{scope}
\path[r] (1) to (7);
\path[r] (7) to (11);
\path[r] (11) to (15);
\path[r] (15) to (1);
\path[r] (1) to (2);
\path[r] (1) to (3);
\path[r] (1) to (4);
\path[r] (7) to (5);
\path[r] (7) to (6);
\path[r] (7) to (8);
\path[r] (11) to (9);
\path[r] (11) to (10);
\path[r] (11) to (12);
\path[r] (15) to (13);
\path[r] (15) to (14);
\path[r] (15) to (16);
\end{tikzpicture}
\end{center}
\caption{The $3$-hairy $4$-cycle $C_{4}\star S_{3}$.}\label{fig:3-hairy}
\end{figure}

In the definition above, the $\star$ notation indicates that we attach a copy of the star $S_m$ at its vertex of degree $m$ to each cycle vertex of $C_n$.  The resulting graph will have $m$ pendants at each cycle vertex. In this case, the tree that we attach to the cycle vertex will be referred to as a \emph{clump}. Figure~\ref{fig:3-hairy} depicts the $3$-hairy $4$-cycle $C_4\star S_3$, which has four clumps, each of which is equal to $S_3$.

Consider the graph $C_n\star S_m$ and let $c_1, c_2, \ldots, c_n$ denote the cycle vertices labeled consecutively.  The general technique for finding a prime vertex labeling for $C_n\star S_m$ is to partition the natural numbers $\mathbb{N}$ into $n$ sets of $m+1$ consecutive integers, denoted $S_1, S_2, \ldots, S_n$. For each $i$ with $1 \leq i \leq n$, assign each $S_i$ to the clump associated to $c_i$. We then hope to find one integer in each set $S_i$ of $n$ consecutive integers that is relatively prime to every other integer in the set and assign it to the cycle vertex $c_i$. In this paper we will focus on $C_n\star S_m$ when $m$ is odd. We begin by tackling $3$-hairy $n$-cycles.

\begin{theorem}\label{thm:labeling 3-hairy}
All $C_{n}\star S_{3}$ are coprime.
\end{theorem}

\begin{proof}
Note that $C_{n}\star S_{3}$ has $4n$ vertices. Let $c_1,c_2,\ldots, c_n$ denote the cycle vertices, and let the three pendant vertices adjacent to $c_i$ be denoted by by $p_{j}^{i}$ for $1 \leq j \leq 3$. Define the labeling function $f:V\to \{1,2,\ldots,4n\}$ via
\begin{align*}
f(c_i) &=\begin{cases}
	1, & i=1\\
    4i-1, & i \geq 2
    \end{cases}\\
f(p^j_i) &=\begin{cases}
	j+1, & i=1,1 \leq j \leq3 \\
	4i-3, & i \geq 2,j=1\\
    4i-2, & i \geq 2,j=2\\
    4i, & i\geq 2,j=3
    \end{cases}.
\end{align*}
Then for $i=1$ and $1 \leq j \leq 3$ it is clear that
\[
(f(c_1),f(p^j_i))=(1,1+j)=1.
\]
For $i \geq 2$, we will show that the remaining pendant vertices have appropriate labels by checking individual values for $j$. If $j=1$, then
\[
(f(c_i),f(p^1_i))=(4i-1,4i-3)=1,
\]
if $j=2$, then
\[
(f(c_i),f(p^2_i))=(4i-1,4i-2)=1,
\]
and lastly if $j=3$, then
\[
(f(c_i),f(p^3_i))=(4i-1,4i)=1.
\]
Finally, to see that all adjacent cycle vertices are assigned relatively prime labels, note that
\[
(f(c_1),f(c_n))=1,
\]
and
\[
(f(c_i),f(c_{i+1}))=(4i-1,4i+3)=1.
\]
Therefore, we can conclude that all 3-hairy $n$-cycles are coprime.
\end{proof}

\begin{figure}[!ht]
\begin{center}
\begin{tikzpicture}[scale=1,auto]
\node (1) at (90:1) [vert] {\scriptsize $1$};
\node (2) at (120:2) [vert] {\scriptsize $2$};
\node (3) at (90:2) [vert] {\scriptsize $3$};
\node (4) at (60:2) [vert] {\scriptsize $4$};
\begin{scope}[rotate=90]
\node (15) at (90:1) [vert] {\scriptsize $15$};
\node (13) at (120:2) [vert] {\scriptsize $13$};
\node (14) at (90:2) [vert] {\scriptsize $14$};
\node (16) at (60:2) [vert] {\scriptsize $16$};
\end{scope}
\begin{scope}[rotate=-90]
\node (7) at (90:1) [vert] {\scriptsize $7$};
\node (5) at (120:2) [vert] {\scriptsize $5$};
\node (6) at (90:2) [vert] {\scriptsize $6$};
\node (8) at (60:2) [vert] {\scriptsize $8$};
\end{scope}
\begin{scope} [rotate=180]
\node (11) at (90:1) [vert] {\scriptsize $11$};
\node (9) at (120:2) [vert] {\scriptsize $9$};
\node (10) at (90:2) [vert] {\scriptsize $10$};
\node (12) at (60:2) [vert] {\scriptsize $12$};
\end{scope}
\path[r] (1) to (7);
\path[r] (7) to (11);
\path[r] (11) to (15);
\path[r] (15) to (1);
\path[r] (1) to (2);
\path[r] (1) to (3);
\path[r] (1) to (4);
\path[r] (7) to (5);
\path[r] (7) to (6);
\path[r] (7) to (8);
\path[r] (11) to (9);
\path[r] (11) to (10);
\path[r] (11) to (12);
\path[r] (15) to (13);
\path[r] (15) to (14);
\path[r] (15) to (16);
\end{tikzpicture}
\end{center}
\caption{An example of a prime labeling for $C_{4}\star S_{3}$.}\label{fig:labeled 3-hairy}
\end{figure}
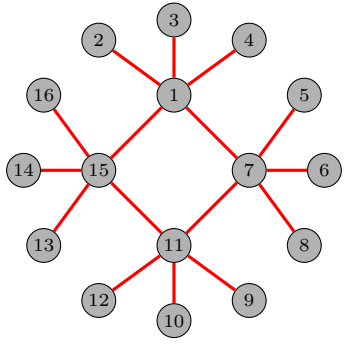

Figure~\ref{fig:labeled 3-hairy} shows the prime vertex labeling for $C_4\star S_3$ that agrees with the labeling described in the proof of Theorem~\ref{thm:labeling 3-hairy}. Next, we address $5$-hairy $n$-cycles.

\begin{theorem}\label{thm:labeling 5-hairy}
All $C_{n}\star S_{5}$ are coprime.
\end{theorem}

\begin{proof}
Note that $C_{n}\star S_{5}$ has $6n$ vertices. Let $c_1, c_2, \ldots, c_n$ denote the cycle vertices, and let the pendant vertices adjacent to to $c_i$ be denoted by by $p_i^j$ for $1\leq j\leq 5$.
Define the labeling function $f: V\to \{1,2,\ldots,6n\}$ via
\begin{align*}
   f(c_i) =& \begin{cases}
       1, & i = 1\\
       6(i-1)+5, & i \geq 2
     \end{cases}\\
   f(p_i^j) =& \begin{cases}
       j+1, & i = 1 , 1\leq j\leq 5\\
       6(i-1)+j, &  i \geq 2 , 1\leq j\leq 4\\
       6(i-1)+6, & i\geq 2 , j=5
     \end{cases}.
\end{align*}
Then for $i=1$ and $1\leq j \leq 5$, it is clear that
\[
(f(c_1),f(p_1^j))=(1,j+1)=1.
\]
Similarly, for $2\leq i \leq n$ and $1\leq j \leq 4$, we have
\[
(f(c_i),f(p_i^j))=(6(i-1)+5, 6(i-1)+j)=1,
\]
and for $2\leq i \leq n$ and $j=5$,
\[
(f(c_i),f(p_i^5))=(6(i-1)+5,6(i-1)+6)=1.
\]
Finally, to see that all adjacent cycle vertices are assigned relatively prime labels, note that
\[
(f(c_1),f(c_2))=(1,6i-1)=1,
\]
and
\[
(f(c_1),f(c_n))=(1,6n-1)=1.
\]
This implies that for $i \geq 2$, we have
\[
(f(c_i),f(c_{i+1}))=(6(i-1)+5,  6(i+1-1)+5)=(6i-1,6i+5)=1.
\]
We have checked all possible adjacencies, and hence all $5$-hairy $n$-cycles are coprime.
\end{proof}

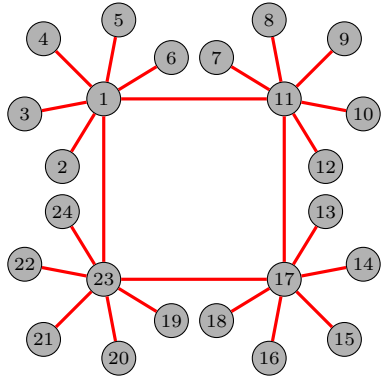
\begin{figure}[!ht]
\begin{center}
\begin{tikzpicture}[scale=1,auto]
\node (1) at (-1.2,1.2) [vert] {\scriptsize $1$};
\node (4) at (-2,2) [vert] {\scriptsize $4$};
\node (5) at (-1,2.25) [vert] {\scriptsize $5$};
\node (6) at (-.3,1.75) [vert] {\scriptsize $6$};
\node (2) at (-1.75,.3) [vert] {\scriptsize $2$};
\node (3) at (-2.25,1) [vert] {\scriptsize $3$};
\path[r] (1) to (2);
\path[r] (1) to (3);
\path[r] (1) to (4);
\path[r] (1) to (5);
\path[r] (1) to (6);
\begin{scope}[rotate=90]
\node (23) at (-1.2,1.2) [vert] {\scriptsize $23$};
\node (21) at (-2,2) [vert] {\scriptsize $21$};
\node (22) at (-1,2.25) [vert] {\scriptsize $22$};
\node (24) at (-.3,1.75) [vert] {\scriptsize $24$};
\node (19) at (-1.75,.3) [vert] {\scriptsize $19$};
\node (20) at (-2.25,1) [vert] {\scriptsize $20$};
\path[r] (23) to (21);
\path[r] (23) to (22);
\path[r] (23) to (24);
\path[r] (23) to (19);
\path[r] (23) to (20);
\end{scope}
\begin{scope}[rotate=180]
\node (17) at (-1.2,1.2) [vert] {\scriptsize $17$};
\node (15) at (-2,2) [vert] {\scriptsize $15$};
\node (16) at (-1,2.25) [vert] {\scriptsize $16$};
\node (18) at (-.3,1.75) [vert] {\scriptsize $18$};
\node (13) at (-1.75,.3) [vert] {\scriptsize $13$};
\node (14) at (-2.25,1) [vert] {\scriptsize $14$};
\path[r] (17) to (15);
\path[r] (17) to (16);
\path[r] (17) to (18);
\path[r] (17) to (13);
\path[r] (17) to (14);
\end{scope}
\begin{scope}[rotate=-90]
\node (11) at (-1.2,1.2) [vert] {\scriptsize $11$};
\node (9) at (-2,2) [vert] {\scriptsize $9$};
\node (10) at (-1,2.25) [vert] {\scriptsize $10$};
\node (12) at (-.3,1.75) [vert] {\scriptsize $12$};
\node (7) at (-1.75,.3) [vert] {\scriptsize $7$};
\node (8) at (-2.25,1) [vert] {\scriptsize $8$};
\path[r] (11) to (7);
\path[r] (11) to (8);
\path[r] (11) to (9);
\path[r] (11) to (10);
\path[r] (11) to (12);
\end{scope}
\path[r] (1) to (23);
\path[r] (23) to (17);
\path[r] (17) to (11);
\path[r] (11) to (1);
\end{tikzpicture}
\end{center}
\caption{An example of a prime labeling for $C_4\star S_5$.}\label{fig:labeled 5-hairy}
\end{figure}

Figure~\ref{fig:labeled 5-hairy} shows the prime vertex labeling for $C_4\star S_5$ that agrees with the labeling described in the proof of Theorem~\ref{thm:labeling 5-hairy}. Continuing with an odd number of pendants, we now handle $C_n\star S_7$.

\begin{theorem}\label{thm:labeling 7-hairy}
All $C_{n}\star S_{7}$ are coprime.
\end{theorem}

\begin{proof}
Note that $C_{n}\star S_{7}$ has $8n$ vertices. Let $c_1,c_2,\ldots, c_n$ denote the cycle vertices, and let the pendant vertices adjacent to $c_i$ be denoted by $p_{i}^{j}$ for $1 \leq j \leq 7$. Define the labeling function $f:V\to \{1,2,\ldots, 8n\}$ via
\begin{align*}
f(c_1)&=1\\
f(p_{1}^{j})&=j+1\\
f(c_i)&=\begin{cases}
	8i-5, & i \equiv_{15} 2,3,6,8,9,11,12,14\\
    8i-3, & i \equiv_{15} 4,5,7,10,13\\
    8i-1, & i \equiv_{15} 0,1
\end{cases}\\
f(p_{i}^{j})&\in \{8i-7,8i-6, \ldots, 8i\}\setminus \{f(c_i)\},
\end{align*}
where each $f(p_{i}^{j})$ is a unique element of $\{8i-7,8i-6, \ldots, 8i\}\setminus \{f(c_i)\}$, the choice being immaterial. In essence, the labeling function splits up $\mathbb{N}$ into smaller sets of eight consecutive numbers each, such as $\{1,2,3,4,5,6,7,8\}$ and $\{9,10,11,12,13,14,15,16\}$. In each set of 8 numbers, the second, third, or fourth odd number is assigned to the cycle vertex such that no multiples of 3 or 5 are chosen, and the rest of the numbers are left for the pendants. The first set is assigned as labels to the first cycle vertex and its associated pendants, the second set is assigned to the second cycle vertex and its associated pendants and so on. This way, one of the numbers from each set is assigned to the cycle vertex, and the rest can be assigned to the pendant vertices in no particular order.

This leaves seven possible cases for which cycle vertices can be adjacent:
\begin{itemize}
\item[] \textbf{Case 1.} $(f(c_i),f(c_{i+1}))=(8i-5,8(i+1)-5)=(8i-5,8i+3)=1$,
\item[] \textbf{Case 2.} $(f(c_i),f(c_{i+1}))=(8i-5,8(i+1)-3)=(8i-5,8i+5)=1$,
\item[] \textbf{Case 3.} $(f(c_i),f(c_{i+1}))=(8i-5,8(i+1)-1)=(8i-5,8i+7)=1$,
\item[] \textbf{Case 4.} $(f(c_i),f(c_{i+1}))=(8i-3,8(i+1)-5)=(8i-3,8i+3)=1$,
\item[] \textbf{Case 5.} $(f(c_i),f(c_{i+1}))=(8i-3,8(i+1)-3)=(8i-3,8i+5)=1$,
\item[] \textbf{Case 6.} $(f(c_i),f(c_{i+1}))=(8i-1,8(i+1)-5)=(8i-1,8i+3)=1$,
\item[] \textbf{Case 7.} $(f(c_i),f(c_{i+1}))=(8i-1,8(i+1)-1)=(8i-1,8i+7)=1$.
\end{itemize}

Cases 1, 5, 6, and 7 correspond to labels that are odd numbers separated by a power of 2. These will always be relatively prime, because if the labels share a common divisor, the divisor must also divide the difference of the labels, a power of two, meaning the divisor is even. This cannot be, as both labels are odd. Cases 3 and 4 correspond to labels that are odd numbers separated by 12 and 6 respectively. Similarly, these are relatively prime: any shared divisor must also divide their difference, 12 or 6, and since neither label is a multiple of 3, the shared divisor would need to be even. Case 2 corresponds to labels that are odd numbers separated by 10, which are relatively prime because neither label is a multiple of 5 or 2.

Showing the selected label for the cycle vertex is relatively prime to the other 7 numbers left for the pendant labels requires three cases:
\begin{itemize}
\item[] \textbf{Case 1.} $f(c_i)=8i-5$,
\item[] \textbf{Case 2.} $f(c_i)=8i-3$,
\item[] \textbf{Case 3.} $f(c_i)=8i-1$.
\end{itemize}

In each of the cases, the fact that a common divisor of two numbers divides their differences is used extensively. In Case 1, the differences between $f(c_i)$ and $f(p_i^j)$ is either 1, 2, 3, 4 or 5. Since the label for $c_i$ is not a multiple of any of these numbers, the cycle vertex label is relatively prime to the labels of the pendant vertices. Shown below is Case 1.
\begin{align*}
(f(c_i),f(p^6_i))&=(8i-5,8i-7)=1,\\
(f(c_i),f(p^2_i))&=(8i-5,8i-6)=1,\\
(f(c_i),f(p^3_i))&=(8i-5,8i-4)=1\\
(f(c_i),f(p^4_i))&=(8i-5,8i-3)=1,\\
(f(c_i),f(p^5_i))&=(8i-5,8i-2)=1,\\
(f(c_i),f(p^6_i))&=(8i-5,8i-1)=1,\\
(f(c_i),f(p^7_i))&=(8i-5,8i)=1.
\end{align*}
Cases 2 and 3 follow similarly. We have checked all possible adjacencies.  Therefore, all 7-hairy $n$-cycles are coprime.
\end{proof}

\begin{figure}[!tht]
\begin{center}
\begin{tikzpicture}[scale=1,auto]
\node (1) at (90:1) [vert] {\scriptsize $1$};
\node (2) at (130:2.5) [vert] {\scriptsize $2$};
\node (3) at (116.67:2.5) [vert] {\scriptsize $3$};
\node (4) at (103.33:2.5) [vert] {\scriptsize $4$};
\node (5) at (90:2.5) [vert] {\scriptsize $5$};
\node (6) at (76.67:2.5) [vert] {\scriptsize $6$};
\node (7) at (63.33:2.5) [vert] {\scriptsize $7$};
\node (8) at (50:2.5) [vert] {\scriptsize $8$};
\begin{scope}[rotate=120]
\node (19) at (90:1) [vert] {\scriptsize $19$};
\node (17) at (130:2.5) [vert] {\scriptsize $17$};
\node (18) at (116.67:2.5) [vert] {\scriptsize $18$};
\node (20) at (103.33:2.5) [vert] {\scriptsize $20$};
\node (21) at (90:2.5) [vert] {\scriptsize $21$};
\node (22) at (76.67:2.5) [vert] {\scriptsize $22$};
\node (23) at (63.33:2.5) [vert] {\scriptsize $23$};
\node (24) at (50:2.5) [vert] {\scriptsize $24$};
\end{scope}
\begin{scope}[rotate=-120]
\node (11) at (90:1) [vert] {\scriptsize $11$};
\node (9) at (130:2.5) [vert] {\scriptsize $9$};
\node (10) at (116.67:2.5) [vert] {\scriptsize $10$};
\node (12) at (103.33:2.5) [vert] {\scriptsize $12$};
\node (13) at (90:2.5) [vert] {\scriptsize $13$};
\node (14) at (76.67:2.5) [vert] {\scriptsize $14$};
\node (15) at (63.33:2.5) [vert] {\scriptsize $15$};
\node (16) at (50:2.5) [vert] {\scriptsize $16$};
\end{scope}
\path[r] (1) to (11);
\path[r] (11) to (19);
\path[r] (19) to (1);
\path[r] (1) to (2);
\path[r] (1) to (3);
\path[r] (1) to (4);
\path[r] (1) to (5);
\path[r] (1) to (6);
\path[r] (1) to (7);
\path[r] (1) to (8);
\path[r] (11) to (9);
\path[r] (11) to (10);
\path[r] (11) to (12);
\path[r] (11) to (13);
\path[r] (11) to (14);
\path[r] (11) to (15);
\path[r] (11) to (16);
\path[r] (19) to (17);
\path[r] (19) to (18);
\path[r] (19) to (20);
\path[r] (19) to (21);
\path[r] (19) to (22);
\path[r] (19) to (23);
\path[r] (19) to (24);
\end{tikzpicture}
\end{center}
\caption{An example of a prime labeling for $C_{3}\star S_{7}$.}\label{fig:C3*S7}
\end{figure}
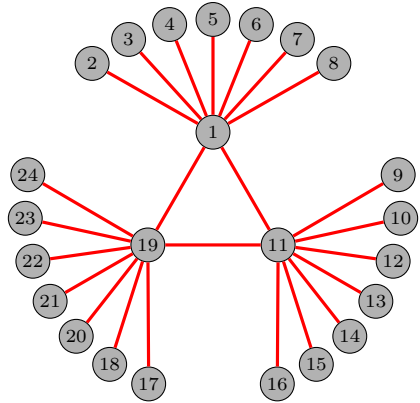

In Figure~\ref{fig:C3*S7}, we have provided a prime vertex labeling for $C_{3}\star S_{7}$ that follows the labeling described in the proof of Theorem~\ref{thm:labeling 7-hairy}.

It is important to note that our work on $m$-hairy $n$-cycles overlaps with the work of Tout et al., who showed, using an existence proof, that all $m$-hairy $n$-cycles are coprime~\cite{Tout1982}. However, our proofs exhibit an explicit prime vertex labeling for $m$ is $3$, $5$, or $7$.

Some families of graphs have natural generalizations based on uniform visual symmetry, such as pendant graphs to $m$-hairy $n$-cycles. But some non-uniform generalizations can be formulated. Consider, for instance, the following number theoretic result.

\begin{proposition}[Bertrand's Postulate]
For every $n\geq 2$, there exists a prime $p$ such that $n < p < 2n$.
\end{proposition}

Using Bertrand's Postulate, we can define the following type of graph, which will naturally yield a prime vertex labeling.

\begin{definition}
Let $n\geq 3$ and consider the cycle $C_n$, where the cycle vertices are consecutively denoted by $c_1, c_2, \ldots, c_n$.  We define the \emph{Bertrand Weed graph}, denoted $BW_n$, to be the non-uniform hairy graph obtained by attaching $2^i - 1$ pendants to each $c_i$.
\end{definition}

\begin{theorem}
All $BW_n$ are coprime.
\end{theorem}

\begin{proof}
By definition of the Bertrand Weed graph, each clump has exactly $2i$ vertices, specifically $2i-1$ from the pendants, and $1$ from the cycle vertex. This allows the natural numbers to be partitioned into sets of size 2, 4, 8, and so on, such that each set is twice as large as the previous. Using Bertrand's Postulate, there is a prime in each set of integers that is assigned to the vertices of each clump in our graph, which is the label given to the corresponding cycle vertex. Thus, all Bertrand weed graphs are coprime.
\end{proof}

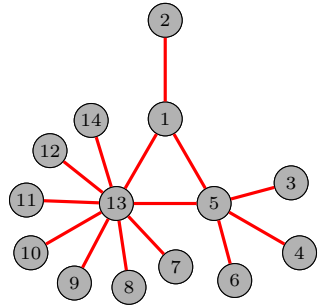
\begin{figure}[!ht]
\begin{center}
\begin{tikzpicture}[scale=.75,auto]
\node (1) at (90:1) [vert] {};
\node (1) at (90:1) [vert] {\scriptsize $1$};
\node (2) at (90:2.75) [vert] {};
\node (2) at (90:2.75) [vert] {\scriptsize $2$};
\path[r] (1) to (2);
\begin{scope}[rotate=120]
\node (13) at (90:1) [vert] {};
\node (13) at (90:1) [vert] {\scriptsize $13$};
\node (10) at (90:2.75) [vert] {};
\node (10) at (90:2.75) [vert] {\scriptsize $10$};
\node (9) at (-.85,2.35) [vert] {};
\node (9) at (-.85,2.35) [vert] {\scriptsize $9$};
\node (11) at (.85,2.35) [vert] {};
\node (11) at (.85,2.35) [vert] {\scriptsize $11$};
\node (14) at (1.5,.65) [vert] {};
\node (14) at (1.5,.65) [vert] {\scriptsize $14$};
\node (7) at (-1.5,.65) [vert] {};
\node (7) at (-1.5,.65) [vert] {\scriptsize $7$};
\node (12) at (1.4,1.55) [vert] {};
\node (12) at (1.4,1.55) [vert] {\scriptsize $12$};
\node (8) at (-1.4,1.55) [vert] {};
\node (8) at (-1.4,1.55) [vert] {\scriptsize $8$};
\path[r] (13) to (10);
\path[r] (13) to (11);
\path[r] (13) to (9);
\path[r] (13) to (14);
\path[r] (13) to (8);
\path[r] (13) to (12);
\path[r] (13) to (7);
\end{scope}
\begin{scope}[rotate=-120]
\node (5) at (90:1) [vert] {};
\node (5) at (90:1) [vert] {\scriptsize $5$};
\node (4) at (90:2.75) [vert] {};
\node (4) at (90:2.75) [vert] {\scriptsize $4$};
\node (3) at (-1,2) [vert] {};
\node (3) at (-1,2) [vert] {\scriptsize $3$};
\node (6) at (1,2) [vert] {};
\node (6) at (1,2) [vert] {\scriptsize $6$};
\path[r] (5) to (4);
\path[r] (5) to (3);
\path[r] (5) to (6);
\end{scope}
\path[r] (1) to (13);
\path[r] (13) to (5);
\path[r] (5) to (1);
\end{tikzpicture}
\end{center}
\caption{An example of a prime labeling for the Bertrand Weed graph $BW_3$.}\label{fig:BW}
\end{figure}

In Figure~\ref{fig:BW}, we have provided a prime vertex labeling for $BW_3$.

\section{Cycles with Complete Ternary Trees}\label{sec:ternary}

This section draws inspiration from the cycles with attached complete binary trees. However, for this family of graphs, cycles will have attached ternary trees instead.

\begin{definition}
A \emph{complete binary tree} is a directed rooted tree with every internal vertex having two children.
\end{definition}

In~\cite{Seoud1999}, Seoud and Youssef showed that every cycle with identical complete binary trees attached to each cycle vertex is coprime. It is natural to ask if cycles with identical complete ternary trees attached to each cycle vertex have prime vertex labelings. For the simplest cases the answer is yes.

\begin{definition}
A \emph{complete ternary tree} is a directed rooted tree with every internal vertex having three children.
\end{definition}

The graphs in Figure~\ref{fig:ternarytrees} depict one and two-level complete ternary trees, respectively. Observe that a one-level complete ternary tree is equal to the star $S_3$.

\tikzstyle{vert} = [circle, draw, fill=grey,inner sep=0pt, minimum size=4.5mm]
\tikzstyle{r} = [draw, very thick, red,-]

\begin{figure}[!ht]
\begin{center}
\begin{tikzpicture}[scale=1,auto]
\node (a) at (0:0) [vert] {\scriptsize $a$};
\node (b) at (180:1) [vert] {\scriptsize $b$};
\node (c) at (90:1) [vert] {\scriptsize $c$};
\node (d) at (0:1) [vert] {\scriptsize $d$};
\path[r] (a) to (b);
\path[r] (a) to (c);
\path[r] (a) to (d);
\begin{scope}[xshift=5cm]
\node (a) at (0:0) [vert] {\scriptsize $a$};
\node (b) at (180:1) [vert] {\scriptsize $b$};
\node (c) at (90:1) [vert] {\scriptsize $c$};
\node (d) at (0:1) [vert] {\scriptsize $d$};
\node (b1) at (180:2) [vert] {\scriptsize $b_1$};
\node (b2) at (160:2) [vert] {\scriptsize $b_2$};
\node (b3) at (200:2) [vert] {\scriptsize $b_3$};
\node (c1) at (90:2) [vert] {\scriptsize $c_1$};
\node (c2) at (70:2) [vert] {\scriptsize $c_2$};
\node (c3) at (110:2) [vert] {\scriptsize $c_3$};
\node (d1) at (0:2) [vert] {\scriptsize $d_1$};
\node (d2) at (-20:2) [vert] {\scriptsize $d_2$};
\node (d3) at (20:2) [vert] {\scriptsize $d_3$};
\path[r] (a) to (b);
\path[r] (a) to (c);
\path[r] (a) to (d);
\path[r] (b) to (b1);
\path[r] (b) to (b2);
\path[r] (b) to (b3);
\path[r] (c) to (c1);
\path[r] (c) to (c2);
\path[r] (c) to (c3);
\path[r] (d) to (d1);
\path[r] (d) to (d2);
\path[r] (d) to (d3);
\end{scope}
\end{tikzpicture}
\end{center}
\caption{Examples of one and two-level complete ternary trees}\label{fig:ternarytrees}
\end{figure}
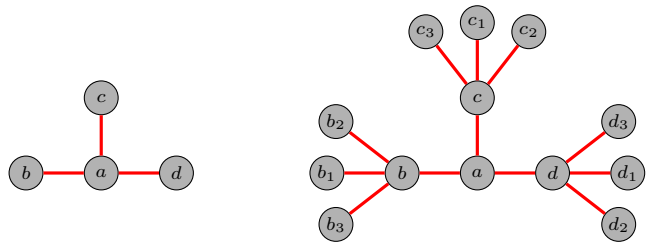

\begin{definition}
We define an \emph{$n$-cycle-pendant with 1-level ternary tree}, denoted $C_n\star P_2\star S_3$, to be the graph that results from first attaching a single pendant to each cycle vertex of $C_n$ followed by attaching a one-level complete ternary tree to each pendant vertex. By extension, we define an \emph{$n$-cycle-pendant with 2-level ternary tree}, denoted $C_n\star P_2\star S_3\star S_3$, to be the graph that results from attaching another complete one-level ternary tree to each vertex of degree one in $C_n\star P_2\star S_3$.
\end{definition}

The star notation here means to first build a cycle, $C_n$, with $n$ vertices. Next, attach one more vertex to each cycle vertex, $c_i$. Now, we have each $c_i$ as one side of a path and one vertex of degree one, which we will denote $p_i$. Finally, attach three more vertices of degree one to each $p_i$, making each $p_i$ the center vertex of a star, $S_3$.

\begin{theorem}\label{thm:labeling 1-level ternary tree}
All $C_n\star P_2\star S_3$ are coprime.
\end{theorem}

\begin{proof}
Note that $C_n\star P_2\star S_3$ contains $5n$ vertices. We will identify our vertices as follows. Let $c_i,1 \leq i \leq n$ denote the cycle vertices, let $p_i$ denote the pendant vertex adjacent to $c_i$, and let the vertices adjacent to $p_i$ be denoted by $s_{i,j}$ for $1 \leq j \leq 3$. Define the labeling function $f:V\to \{1,2, \ldots, 5n\}$ via
\begin{align*}
f(c_i)&= 5i - 4, 1 \leq i \leq n\\
f(p_i)&=\begin{cases}
            5i-2, & i$ is odd $\\
            5i-3, & i\equiv_{2}0,i\not\equiv_{6}0\\
            5i-1, & i\equiv_{6}0
            \end{cases}\\
f(s_{i,j})&=\begin{cases}
            5i - 3 + j, & i \text{ is even}\\
            5i - 2 + j, & j \neq 3 \text{ and } i \text{ is odd} \\
            5i - 3, & j = 3 \text{ and } i \text{ is odd}
            \end{cases}.
\end{align*}
It is straightforward to show that this mapping is injective and all adjacent vertices have relatively prime labels.
\end{proof}

Figure~\ref{fig:labeled 1-level ternary tree} shows the prime vertex labeling for $C_3\star P_2\star S_3$ that agrees with the labeling described in the proof of Theorem~\ref{thm:labeling 1-level ternary tree}.

\tikzstyle{vert} = [circle, draw, fill=grey,inner sep=0pt, minimum size=4.5mm]
\tikzstyle{b} = [draw,very thick,blue,-]
\tikzstyle{r} = [draw, very thick, red,-]
\tikzstyle{g} = [draw, very thick, green, -]
\tikzstyle{blk} = [draw, very thick, black, -]

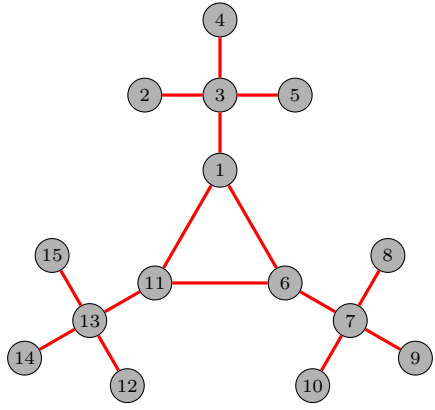
\begin{figure}
\begin{center}
\begin{tikzpicture}[scale=1,auto]
\node (a) at (90:1) [vert] {\scriptsize $1$};
\node (b) at (90:2) [vert] {\scriptsize $3$};
\node (c) at (90:3) [vert] {\scriptsize $4$};
\node (d) at (-1,2) [vert] {\scriptsize $2$};
\node (e) at (1,2) [vert] {\scriptsize $5$};
\path[r] (a) to (b);
\path[r] (b) to (c);
\path[r] (b) to (e);
\path[r] (b) to (d);
\begin{scope}[rotate=120]
\node (a') at (90:1) [vert] {\scriptsize $11$};
\node (b') at (90:2) [vert] {\scriptsize $13$};
\node (c') at (90:3) [vert] {\scriptsize $14$};
\node (d') at (-1,2) [vert] {\scriptsize $12$};
\node (e') at (1,2) [vert] {\scriptsize $15$};
\path[r] (a') to (b');
\path[r] (b') to (c');
\path[r] (b') to (e');
\path[r] (b') to (d');
\end{scope}
\begin{scope}[rotate=-120]
\node (a'') at (90:1) [vert] {\scriptsize $6$};
\node (b'') at (90:2) [vert] {\scriptsize $7$};
\node (c'') at (90:3) [vert] {\scriptsize $9$};
\node (d'') at (-1,2) [vert] {\scriptsize $8$};
\node (e'') at (1,2) [vert] {\scriptsize $10$};
\path[r] (a'') to (b'');
\path[r] (b'') to (c'');
\path[r] (b'') to (e'');
\path[r] (b'') to (d'');
\end{scope}
\path[r] (a) to (a');
\path[r] (a') to (a'');
\path[r] (a'') to (a);
\end{tikzpicture}
\end{center}
\caption{An example of a prime labeling for $C_3\star P_2\star S_3$.}\label{fig:labeled 1-level ternary tree}
\end{figure}

Recall from above that an $n$-cycle-pendant with 2-level ternary tree, denoted $C_n\star P_2\star S_3\star S_3$, is the graph that results from gluing a copy of $S_3$, at the center vertex, onto each of the vertices of degree $1$ in the graph $C_n\star P_2\star S_3$.

\begin{theorem}\label{thm:labeling 2-level ternary tree}
All $C_n\star P_2\star S_3 \star S_3$ are coprime.
\end{theorem}

\begin{proof}
Note that $C_n\star P_2\star S_3 \star S_3$ contains $14n$ vertices. We will identify our vertices as follows. Let $c_i,1 \leq i \leq n$ denote the cycle vertices, let $p_i$ denote the pendant vertex adjacent to $c_i$, let the non-cycle vertices adjacent to $p_i$ be denoted $s_{i,j}$ for $1 \leq j \leq 3$, and let the remaining vertices adjacent to each $s_{i,j}$ be denoted $l_{i,j,k}$ for $1 \leq k \leq 3$. Our labeling function $f:V\to \{1,2, \ldots, 14n\}$ is best defined by first describing cycle and pendent vertex labels:
\begin{align*}
f(c_i)&= 14i - 13, 1 \leq i \leq n\\
f(p_i)&=\begin{cases}
            14i-12, & i \equiv_{3}1,2\\
            14i-10, & i \equiv_{3}0
            \end{cases}.
\end{align*}
The remaining vertex labels are determined by the values of $i$, $j$ and $k$ as follows. If $i\equiv_{3}1,2$, then define
\begin{align*}
f(s_{i,j})&=\begin{cases}
            14i - 9, & j = 1\\
            14i - 5, & j = 2\\
            14i - 3, & j = 3
            \end{cases}\\
f(l_{i,j,k})&=\begin{cases}
            14i - 11, & j = 1, k = 1\\
            14i - 10, & j = 1, k = 2\\
            14i - 8, & j = 1, k = 3\\
            14i - 7, & j = 2, k = 1\\
            14i - 6, & j = 2,k = 2\\
            14i - 4, & j = 2, k = 3\\
            14i - 2, & j = 3, k = 1\\
            14i - 1, & j = 3, k = 2\\
            14i, & j = 3, k = 3
            \end{cases}.
\end{align*}
If $i \equiv_{3}0$, then define
\begin{align*}
f(s_{i,j})&=\begin{cases}
            14i - 11, & j = 1\\
            14i - 7, & j = 2\\
            14i - 1, & j = 3
            \end{cases}\\
f(l_{i,j,k})&=\begin{cases}
            14i - 12, & j = 1, k = 1\\
            14i - 9, & j = 1, k = 2\\
            14i - 8, & j = 1, k = 3\\
            14i - 6, & j = 2, k = 1\\
            14i - 5, & j = 2, k = 2\\
            14i - 4, & j = 2, k = 3\\
            14i - 3, & j = 3, k = 1\\
            14i - 2, & j = 3, k = 2\\
            14i, & j = 3, k = 3
            \end{cases}.
\end{align*}
Again, it is relatively straightforward to check that this mapping results in a prime vertex labeling. The details are left to the interested reader.
\end{proof}

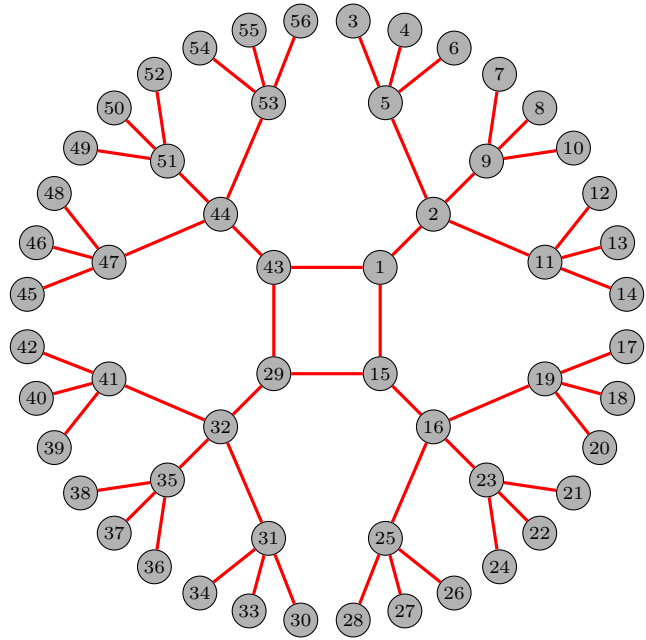
\begin{figure}[!ht]
\begin{center}
\begin{tikzpicture}[scale=1,auto]
\node (a) at (45:1) [vert] {\scriptsize $1$};
\node (b) at (45:2) [vert] {\scriptsize $2$};
\node (c) at (75:3) [vert] {\scriptsize $5$};
\node (d) at (45:3) [vert] {\scriptsize $9$};
\node (e) at (15:3) [vert] {\scriptsize $11$};
\node (f) at (85:4) [vert] {\scriptsize $3$};
\node (g) at (75:4) [vert] {\scriptsize $4$};
\node (h) at (65:4) [vert] {\scriptsize $6$};
\node (i) at (55:4) [vert] {\scriptsize $7$};
\node (j) at (45:4) [vert] {\scriptsize $8$};
\node (k) at (35:4) [vert] {\scriptsize $10$};
\node (l) at (25:4) [vert] {\scriptsize $12$};
\node (m) at (15:4) [vert] {\scriptsize $13$};
\node (n) at (5:4) [vert] {\scriptsize $14$};
\path[r] (a) to (b);
\path[r] (b) to (c);
\path[r] (b) to (e);
\path[r] (b) to (d);
\path[r] (c) to (f);
\path[r] (c) to (g);
\path[r] (c) to (h);
\path[r] (d) to (i);
\path[r] (d) to (j);
\path[r] (d) to (k);
\path[r] (e) to (l);
\path[r] (e) to (m);
\path[r] (e) to (n);
\begin{scope}[rotate=-90]
\node (a') at (45:1) [vert] {\scriptsize $15$};
\node (b') at (45:2) [vert] {\scriptsize $16$};
\node (c') at (75:3) [vert] {\scriptsize $19$};
\node (d') at (45:3) [vert] {\scriptsize $23$};
\node (e') at (15:3) [vert] {\scriptsize $25$};
\node (f') at (85:4) [vert] {\scriptsize $17$};
\node (g') at (75:4) [vert] {\scriptsize $18$};
\node (h') at (65:4) [vert] {\scriptsize $20$};
\node (i') at (55:4) [vert] {\scriptsize $21$};
\node (j') at (45:4) [vert] {\scriptsize $22$};
\node (k') at (35:4) [vert] {\scriptsize $24$};
\node (l') at (25:4) [vert] {\scriptsize $26$};
\node (m') at (15:4) [vert] {\scriptsize $27$};
\node (n') at (5:4) [vert] {\scriptsize $28$};
\path[r] (a') to (b');
\path[r] (b') to (c');
\path[r] (b') to (e');
\path[r] (b') to (d');
\path[r] (c') to (f');
\path[r] (c') to (g');
\path[r] (c') to (h');
\path[r] (d') to (i');
\path[r] (d') to (j');
\path[r] (d') to (k');
\path[r] (e') to (l');
\path[r] (e') to (m');
\path[r] (e') to (n');
\end{scope}
\begin{scope}[rotate=180]
\node (a'') at (45:1) [vert] {\scriptsize $29$};
\node (b'') at (45:2) [vert] {\scriptsize $32$};
\node (c'') at (75:3) [vert] {\scriptsize $31$};
\node (d'') at (45:3) [vert] {\scriptsize $35$};
\node (e'') at (15:3) [vert] {\scriptsize $41$};
\node (f'') at (85:4) [vert] {\scriptsize $30$};
\node (g'') at (75:4) [vert] {\scriptsize $33$};
\node (h'') at (65:4) [vert] {\scriptsize $34$};
\node (i'') at (55:4) [vert] {\scriptsize $36$};
\node (j'') at (45:4) [vert] {\scriptsize $37$};
\node (k'') at (35:4) [vert] {\scriptsize $38$};
\node (l'') at (25:4) [vert] {\scriptsize $39$};
\node (m'') at (15:4) [vert] {\scriptsize $40$};
\node (n'') at (5:4) [vert] {\scriptsize $42$};
\path[r] (a'') to (b'');
\path[r] (b'') to (c'');
\path[r] (b'') to (e'');
\path[r] (b'') to (d'');
\path[r] (c'') to (f'');
\path[r] (c'') to (g'');
\path[r] (c'') to (h'');
\path[r] (d'') to (i'');
\path[r] (d'') to (j'');
\path[r] (d'') to (k'');
\path[r] (e'') to (l'');
\path[r] (e'') to (m'');
\path[r] (e'') to (n'');
\end{scope}
\begin{scope}[rotate=90]
\node (a''') at (45:1) [vert] {\scriptsize $43$};
\node (b''') at (45:2) [vert] {\scriptsize $44$};
\node (c''') at (75:3) [vert] {\scriptsize $47$};
\node (d''') at (45:3) [vert] {\scriptsize $51$};
\node (e''') at (15:3) [vert] {\scriptsize $53$};
\node (f''') at (85:4) [vert] {\scriptsize $45$};
\node (g''') at (75:4) [vert] {\scriptsize $46$};
\node (h''') at (65:4) [vert] {\scriptsize $48$};
\node (i''') at (55:4) [vert] {\scriptsize $49$};
\node (j''') at (45:4) [vert] {\scriptsize $50$};
\node (k''') at (35:4) [vert] {\scriptsize $52$};
\node (l''') at (25:4) [vert] {\scriptsize $54$};
\node (m''') at (15:4) [vert] {\scriptsize $55$};
\node (n''') at (5:4) [vert] {\scriptsize $56$};
\path[r] (a''') to (b''');
\path[r] (b''') to (c''');
\path[r] (b''') to (e''');
\path[r] (b''') to (d''');
\path[r] (c''') to (f''');
\path[r] (c''') to (g''');
\path[r] (c''') to (h''');
\path[r] (d''') to (i''');
\path[r] (d''') to (j''');
\path[r] (d''') to (k''');
\path[r] (e''') to (l''');
\path[r] (e''') to (m''');
\path[r] (e''') to (n''');
\end{scope}
\path[r] (a) to (a');
\path[r] (a) to (a''');
\path[r] (a'') to (a');
\path[r] (a'') to (a''');
\end{tikzpicture}
\end{center}
\caption{An example of a prime labeling for $C_4\star P_2\star S_3\star S_3$.}\label{fig:labeled 2-level ternary tree}
\end{figure}

Figure~\ref{fig:labeled 2-level ternary tree} shows the prime vertex labeling for $C_4\star P_2\star S_3\star S_3$ that agrees with the labeling described in the proof of Theorem~\ref{thm:labeling 2-level ternary tree}.

\section{Conclusion}\label{sec:conclusion}

Seoud and Youssef's conjecture is still open. This makes an excellent target for further work, as well as attempting to find more families of coprime graphs. Both hairy cycles and cycle pendant stars still have more progress that can be made. Specifically, we conjecture that a similar approach will work for labeling hairy cycles up to fifteen pendants and cycle pendant stars up to fourteen outer vertices.

The reasoning for labeling techniques for hairy cycles failing at sixteen ``hairs" stems from the following result proved by Pillai~\cite{Pillai1941}.

\begin{proposition}
When $m \geq 17$, we can find $m$ consecutive integers such that no number in the set is prime to all the rest in the set.
\end{proposition}

That is, according to Pillai, in any set of $17$ or more consecutive integers, there is a subset of those integers in which no element is relatively prime to the rest of the elements in the set. This would imply that at some point, the proposed labeling for hairy cycles would fail. Thus, a new labeling scheme would need to be devised to label hairy cycles when the number of pendants becomes large.

Naturally, the progress we have achieved will not prove Seoud and Youssef's conjecture. Doing so would require a drastically different approach, namely not trying to show specific families of graphs to be coprime. However, in~\cite{Seoud2012}, Seoud et al.~detail necessary and sufficient conditions for a graph to be coprime. For readers more interested in graph labeling in general, additional information can be found in Gallian's dynamic survey on graph labelings~\cite{Gallian2014}.

\section*{Acknowledgements}
We would like to thank our research mentors Dana C.~Ernst and Jeff Rushall for their guidance on this project.

\bibliographystyle{amsplain}
\bibliography{PrimeLabelings}

\end{document}